\documentclass[12pt]{article}
\usepackage[colorlinks=true,
linkcolor=webgreen,
filecolor=webbrown,
citecolor=webgreen]{hyperref}
\usepackage{psfig}

\definecolor{webgreen}{rgb}{0,.5,0}
\definecolor{webbrown}{rgb}{.6,0,0}

\usepackage{amssymb,amsmath,latexsym,psfig,epsf}
\usepackage{showlabels}

\newtheorem{theorem}{Theorem}[section]

\newtheorem{coro}[theorem]{Corollary}

\newcommand{\be}{\begin{enumerate}}
\newcommand{\ee}{\end{enumerate}}
\newcommand{\bi}{\begin{itemize}}
\newcommand{\ei}{\end{itemize}}

\newcommand{\ba}{\begin{array}}
\newcommand{\ea}{\end{array}}

\newcommand{\Z}{{\mathbb{Z}}}

\newcommand{\RR}{{\mathbb{R}}}
\newcommand{\FF}{{\mathbb{F}}}
\newcommand{\CL}{{\cal{C}}}
\newcommand{\CC}{\mathbb C}
\renewcommand{\em}{\sf}

\makeatletter
\def\@sect#1#2#3#4#5#6[#7]#8{\ifnum #2>\c@secnumdepth
     \def\@svsec{}\else 
     \refstepcounter{#1}\edef\@svsec{\csname the#1\endcsname.\hskip .75em }\fi
     \@tempskipa #5\relax
      \ifdim \@tempskipa>\z@ 
        \begingroup #6\relax
          \@hangfrom{\hskip #3\relax\@svsec}{\interlinepenalty \@M #8\par}%
        \endgroup
       \csname #1mark\endcsname{#7}\addcontentsline
         {toc}{#1}{\ifnum #2>\c@secnumdepth \else
                      \protect\numberline{\csname the#1\endcsname}\fi
                    #7}\else
        \def\@svsechd{#6\hskip #3\@svsec #8\csname #1mark\endcsname
                      {#7}\addcontentsline
                           {toc}{#1}{\ifnum #2>\c@secnumdepth \else
                             \protect\numberline{\csname the#1\endcsname}\fi
                       #7}}\fi
     \@xsect{#5}}
\def\@begintheorem#1#2{\it \trivlist \item[\hskip \labelsep{\bf #1\ #2.}]}
\def\section{\@startsection {section}{1}{\z@}{-3.5ex plus -1ex minus 
 -.2ex}{2.3ex plus .2ex}{\normalsize\bf}}
\makeatother

\begin{document}
\begin{center}
{\Large {\bf  A Simple Construction for the Barnes-Wall Lattices}} \\
\vspace{1.5\baselineskip}
{\em Gabriele Nebe} \\
\vspace*{1\baselineskip}
Abteilung Reine Mathematik \\
Universitaet Ulm, 89069 Ulm, Germany \\
\vspace{1\baselineskip}
and \\
\vspace{1\baselineskip}
{\em E. M. Rains} and {\em N. J. A. Sloane} \\
\vspace*{1\baselineskip}
Information Sciences Research, AT\&T Shannon Labs \\
180 Park Avenue, Florham Park, NJ 07932-0971, U.S.A. \\
\vspace{1.5\baselineskip}
To Dave Forney, on the occasion of his sixtieth birthday \\
\vspace{1.5\baselineskip}
{\bf ABSTRACT} \\
\vspace{.5\baselineskip}
\end{center}
\setlength{\baselineskip}{1.5\baselineskip}

A certain family of orthogonal groups (called ``Clifford groups'' by G. E. Wall) has arisen in a variety of different contexts in recent years.
These groups have a simple definition as the automorphism groups of certain
generalized Barnes-Wall lattices.
This leads to an especially simple construction for the usual Barnes-Wall lattices.

\{This is based on the third author's talk at the Forney-Fest, M.I.T., 
March 2000, which in turn is based on our paper
``\htmladdnormallink{The Invariants of the Clifford
Groups}{http://www.research.att.com/~njas/doc/cliff1.pdf}'' 
{\it Designs, Codes, Crypt.}, {\bf 24} (2001), 99--121,
to which the reader is referred for further details and proofs.\}

\vspace{1\baselineskip}

\clearpage
\setcounter{page}{2}

\section{Background}
The Barnes-Wall lattices define
an infinite sequence of sphere packings in dimensions $2^m$, $m \ge 0$, which include the densest packings known in dimensions 1, 2, 4, 8 and 16 \cite{BW59}, \cite{SPLAG}.
In dimensions 32 and higher they are less dense than other known packings, but they are still interesting for other reasons --- they form one of the few infinite sequences of lattices where it is possible to do explicit calculations.
For example, there is an explicit formula for their kissing numbers \cite{SPLAG}.
This talk will describe a beautifully simple construction for these lattices that we found in the summer of 1999.
A more comprehensive account will appear elsewhere \cite{cliff1}, \cite{cliff2}.
Since Dave Forney is fond of the Barnes-Wall lattices (cf. \cite{For88}, \cite{For88a}) we hope he will like this construction as much as we do.

This work had its origin in 1995 when J. H. Conway, R. H. Hardin and N. J. A. S. were studying packings in Grassmann manifolds --- in other words,
packings of Euclidean $k$-dimensional subspaces in $n$-dimensional space \cite{grass1}.
One of our nicest constructions was an optimal packing of 70 4-dimensional subspaces in $\RR^8$.
The symmetry group of this packing (the subgroup of the orthogonal group $O(8, \RR )$ that fixes the collection of subspaces) has order 5160960.

Shortly afterwards, the {\em same} 8-dimensional group arose in the work of P. W. Shor and others on quantum computers (cf. \cite{BDSW96}, \cite{Kit97}).

This astonishing coincidence --- see \cite{QC2} for the full story --- drew
attention to earlier work on the family to which this group belongs
\cite{Wall3}, \cite{Wall1}, \cite{Wall2}, \cite{CCKS}, \cite{Wall4}.
Following Wall, we call these Clifford groups, although these are not the groups usually referred
to by this name \cite{Port}.
Investigation of the representations of subgroups of these groups led to further constructions of optimal packings in Grassmann manifolds \cite{grass3} and constructions of quantum error-correcting codes \cite{QC1}, \cite{QC2}.

Independently, and around the same time, these groups\footnote{Although
at that time they were not recognized as the Clifford groups.} also occurred in the work of V. M. Sidelnikov, in connection with the construction of spherical $t$-designs \cite{Kaz}, \cite{Sid1}, \cite{Sid2}, \cite{Sid3}, \cite{Sid4}.

The complete account of our work \cite{cliff1}, \cite{cliff2} describes the invariants of these Clifford groups and their connections with binary self-dual codes.
Much of this work had been anticipated by Runge \cite{Schottky},
\cite{Rung93}, \cite{Rung95}, \cite{Rung96}.
In our two papers we
clarify the connections with spherical $t$-designs and Sidelnikov's work, and also generalize these results to the complex Clifford groups and doubly-even binary self-dual codes.
Again the main result was first given by Runge.

In recent years many other kinds of self-dual codes have been studied by a number of authors.
Nine such families were named and surveyed in \cite{chapter}.
In \cite{cliff2} we give a general definition of the ``type'' of a self-dual code which includes all these families as well as other self-dual codes over rings and modules.
For each ``type'' we investigate the structure of the associated ``Clifford-Weil group'' and its ring of invariants.

Some of the results in \cite{cliff1}, \cite{cliff2} can be regarded as providing a general setting for Gleason's theorems
\cite{Gleason}, \cite{MS77}, \cite{chapter} about the weight enumerator of a binary self-dual code, a doubly-even binary self-dual code and a self-dual code over $\FF_p$.
They are also a kind of discrete analogue of a long series of theorems going back to Eichler (see for example \cite{Boech}, \cite{Rung93}, \cite{Rung95}, \cite{Rung96}), stating that under certain conditions theta series of quadratic forms are bases for spaces of modular forms:
here complete weight enumerators of generalized self-dual codes are bases for spaces of invariants of ``Clifford-Weil groups''.

\section{A simple construction for the Barnes-Wall lattices}
\begin{figure}[htb]
\caption{The two Barnes-Wall lattices $L_1$ (solid circles) and $L'_1$ (solid or hollow circles) in two dimensions.}

\vspace*{+.2in}
\begin{center}
\epsfxsize=2.5in
\leavevmode\epsffile{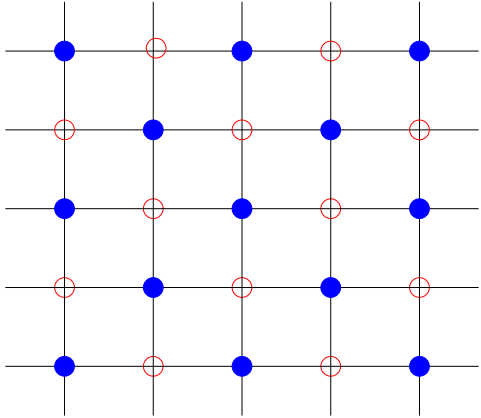}
\end{center}

\label{F1}
\end{figure}

There is a pair of Barnes-Wall lattices $L_m$ and $L'_m$ in each dimension
$2^m$, $m \ge 0$.
The two lattices are geometrically similar\footnote{In other words they differ only by a rotation and change of scale.} and $L_m$
is a sublattice of index  $2^k$, $k = 2^{m-1}$,
in $L'_m$.
In dimension 2 these lattices are shown in Fig. \ref{F1}, where $L_1$ consists
of the points marked with solid circles and $L'_1$ consists of the points marked with either solid or hollow circles.
Both are geometrically similar to the square lattice $\Z^2$.

\begin{figure}[htb]
\caption{The eight minimal vectors of $L_1$ and $\sqrt{2}$ $L'_1$.
The $\Z [\sqrt{2}]$ span of these points is the ``balanced'' Barnes-Wall lattice $M_1$.}

\vspace*{+.2in}
\begin{center}
\epsfxsize=1.5in
\leavevmode\epsffile{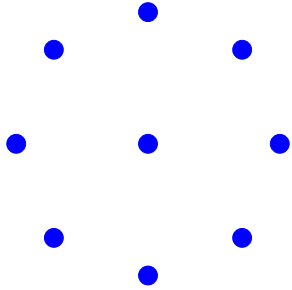}
\end{center}

\label{F2}
\end{figure}

Suppose we multiply the points of $L'_1$ by $\sqrt{2}$.
Then the eight minimal vectors of $L_1$ and $\sqrt{2}$ $L'_1$ now have the same length and form the familiar
configuration of points used in the 8-PSK signaling system (Fig.~\ref{F2}).

We now define the generalized or ``balanced'' Barnes-Wall lattice $M_1$ to be the set of all $\Z [\sqrt{2}]$-integer combinations of the eight vectors in Figure 2.
That is, we take integer combinations of these vectors where ``integer'' now means a number of the form $a+b \sqrt{2}$, $a,b \in \Z$.
In more formal language, $M_1$ is a $\Z[\sqrt{2}]$-lattice (or $\Z [\sqrt{2}]$-module).
Note that we can recover $L_1$ from $M_1$ by taking just those vectors in $M_1$ whose components are integers.

In general we define
the rational part of a $\Z [\sqrt{2}]$-lattice $\Lambda$ to
consist of the vectors which have rational components, and
the irrational part to consist of the vectors whose components are rational 
multiples of $\sqrt{2}$.
We can now state the construction.

\begin{theorem}\label{th1}
Define the balanced Barnes-Wall lattice $M_m$ to be the $m$-fold tensor product $M_1^{\otimes m}$.
Then the rational part of $M_m$ is the
Barnes-Wall lattice $L_m$, and the purely irrational part is $\sqrt{2}$ $L'_m$.
\end{theorem}

For the proof see \cite{cliff1}.

To be quite explicit, note that we need only two of the vectors in Fig.~\ref{F2}, and we can take
$$G_1 = \left( \begin{array}{cc}
\sqrt{2} & 0 \\
1 & 1
\end{array}\right)$$
as a generator matrix for $M_1$.
Then the $m$-fold tensor power of this matrix,
$$G_m = G_1^{\otimes m} = G_1 \otimes G_1 \otimes \cdots \otimes G_1$$
is a generator matrix for $M_m$.

For example, $G_2 = G_1 \otimes G_1$ is
$$
\left[
\begin{array}{cccc}
2 & 0 & 0 & 0 \\
\sqrt{2} & \sqrt{2} & 0 & 0 \\
\sqrt{2} & 0 & \sqrt{2} & 0 \\
1 & 1 & 1 & 1
\end{array}
\right] ~.
$$
The rational part, $L_2$, is generated by
$$
\left[\begin{array}{cccc}
2 & 0 & 0 & 0 \\
2 & 2 & 0 & 0 \\
2 & 0 & 2 & 0 \\
1 & 1 & 1 & 1
\end{array}
\right]\,,
~~\mbox{or equivalently}~~
\left[\begin{array}{cccc}
2 & 0 & 0 & 0 \\
0 & 2 & 0 & 0 \\
0 & 0 & 2 & 0 \\
1 & 1 & 1 & 1
\end{array}
\right] \,.
$$
This lattice is geometrically similar to $D_4$ (\cite{SPLAG}, Chap. 4, Eq. (90)).
The purely irrational part, $\sqrt{2}L'_2$, is generated by
$$
\left[\begin{array}{cccc}
2 \sqrt{2} & 0 & 0 & 0 \\
\sqrt{2} & \sqrt{2} & 0 & 0 \\
\sqrt{2} & 0 & \sqrt{2} & 0 \\
\sqrt{2} & \sqrt{2} & \sqrt{2} & \sqrt{2}
\end{array}
\right]\,,
~~\mbox{or equivalently}~~
\sqrt{2}
\left[\begin{array}{cccc}
2 & 0 & 0 & 0 \\
1 & 1 & 0 & 0 \\
1 & 0 & 1 & 0 \\
1 & 0 & 0 & 1
\end{array}
\right]\,,
$$
which another version of $D_4$ (\cite{SPLAG}, Chap. 4, Eq. (86)).

We may avoid the use of coordinates and work directly with Gram matrices
or quadratic forms, provided we select an appropriate element $\phi$
of the Galois group.
Let $u_1$ and $u_2$ be the generating vectors corresponding to the rows of the matrix $G_1$, and let $\phi$ negate $u_1$ and fix $u_2$.

Then $M_1$ is the $\Z [\sqrt{2}]$-lattice with Gram matrix
$$A_1 = G_1 G_1^{tr} = \left(
\begin{array}{cc}
2 & \sqrt{2} \\
\sqrt{2} & 2
\end{array}
\right) \,,$$
$L_1$ is the sublattice of $M_1$ fixed by $\phi$ and $\sqrt{2}$ $L'_1$ is the sublattice negated by $\phi$.

Furthermore,
$\phi$ has a natural extension to $M_m = M_1^{\otimes m}$, and $L_m$ is the sublattice of $M_n$ fixed by $\phi$ and $\sqrt{2}$ $L'_m$ is the sublattice negated by $\phi$.

\section{The Clifford groups and their invariants}
The Clifford groups $\CL_m$ mentioned at the beginning of the paper now have a very simple definition:
for all $m \ge 1$, $\CL_m$ is $Aut (M_n)$, i.e. the subgroup of $O(2^m, \RR )$ that preserves $M_m$.

For the proof that this definition is equivalent to the usual one given in \cite{CCKS}, \cite{grass3}, \cite{QC2}, see Proposition 5.3 of \cite{cliff1}.

An {\em invariant polynomial} of $\CL_m$ is a polynomial in $2^m$ variables with real coefficients that is fixed by every element of the group \cite{Benson}.
The ring of invariant polynomials plays an important role in constructing spherical $t$-designs from the group (see for example \cite{SPLAG}, Chap. 3, Section 4.2).
The {\em Molien series} of the group is a generating function for the numbers of linearly
independent homogeneous invariants of each degree \cite{Benson}, \cite{MS77}, \cite{chapter}.

In \cite{CCKS} it was asked ``is it possible to say something about the
Molien series [of the groups $\CL_m$],
such as the minimal degree of an invariant?''
Such questions also arise in the work of
Sidelnikov \cite{Sid1}, \cite{Sid2}, \cite{Sid3}, \cite{Sid4}.
The answers are given by the following theorem of Runge \cite{Schottky},
\cite{Rung93}, \cite{Rung95}, \cite{Rung96}.

\begin{theorem}\label{th2}
$($Runge; \cite{cliff1}$)$.
Fix integers $k$ and $m \ge 1$.
The space of homogeneous invariants of $\CL_m$ of degree $2k$ is spanned by the complete weight enumerators of the codes $C \otimes GF (2^m)$, where $C$ ranges over all binary
self-dual codes of length $2k$;
this is a basis if $m \ge k-1$.
\end{theorem}

We rediscovered this result in the summer of 1999.
Our proof is somewhat simpler than Runge's as it avoids the use of Siegel modular forms \cite{cliff1}.

\begin{coro}\label{th3}
Let $\Phi_m (t)$ be the Molien series of $\CL_m$.
As $m$
tends to infinity, the series $\Phi_m(t)$ tend monotonically
to
\[
\sum_{k = 0}^{\infty} N_{2k} t^{2k} ~,
\]
where $N_{2k}$ is the number of equivalence classes of binary self-dual codes of
length $2k$.
\end{coro}

Explicit calculations for $m=1,2$ show:

\begin{coro}\label{th4}
The initial terms of the Molien series of $\CL_m$
are given by
$$
1+t^2+t^4+t^6+2 t^8 + 2 t^{10} + O(t^{12})\,,
$$
where the next term is $2 t^{12}$ for $m=1$, and $3 t^{12}$ for $m>1$.
\end{coro}

Sidelnikov \cite{Sid2}, \cite{Sid3} showed that the lowest degree of a harmonic
invariant
of ${\cal C}_m$ is $8$.
Inspection of the above Molien series gives the following
stronger result.

\begin{coro}{\label{th5}}
The smallest degree of a harmonic invariant of ${\cal C}_m$ is $8$, and
there is a unique harmonic invariant of degree $8$.
There are no harmonic invariants of degree $10$.
\end{coro}

There is a unique harmonic invariant of degree 8, which can be taken to be the complete weight enumerator of $H_8 \otimes GF(2^m)$, where $H_8$ is the $[8,4,4]$ binary Hamming code,
minus a suitable multiple of the fourth power of the quadratic form.

From Corollary \ref{th5} we can easily show that appropriate orbits under $\CL_m$ form spherical 7-designs, 11-designs, etc.

We conclude by mentioning one last result (Corollary 5.7 of \cite{cliff1}).

\begin{theorem}\label{th6}
Let $C$ be any binary self-dual code that is not generated by vectors of weight 2,
and form the complete weight enumerator of $C \otimes GF(2^m)$.
Then the subgroup of $O(2^m, \RR )$ that fixes this weight enumerator is precisely $\CL_m$.
\end{theorem}

As already mentioned, there are analogues of all these results for the complex versions of the Clifford group.
Now ``self-dual code'' is replaced by ``doubly-even self-dual code''.

The case $m=1$ of Theorem \ref{th6} and its complex analogue imply the following.

Let $C$ be a binary self-dual code of length $N$ and $W(x,y)$ its Hamming weight enumerator.
Let $G$ be the subgroup of $O(2, \RR)$ that fixes $W(x,y)$.
Provided $C$ is not generated by vectors of weight 2, $G \cong \CL_1$, of order 16.
If $C$ is doubly-even, the subgroup of the unitary group $U(2, \CC)$ that fixes $W(x,y)$ is (apart from its center, which of course may contain complex $N$-th roots of unity)
the familiar group of order 192 arising in Gleason's theorem.

It was of course known that $W(x,y)$ is fixed by these two groups, of order 16 and 192.
But we were not aware before this of any proof that the group of $W(x,y)$ could never be bigger.

\end{document}